  \documentclass[10pt]{amsart}
  \usepackage{epic,latexsym}

 \newtheorem{theorem}{Theorem}[section]
 
 \newtheorem{definition}[theorem]{Definition}

 \newtheorem{remark}[theorem]{Remark}
 \newtheorem{example}[theorem]{Example}
 \newtheorem{conjecture}[theorem]{Conjecture}

  \newcommand{\remind}[1]{{}}

  \newcommand{\cited}{{}}

\newcommand{\bpf}{\noindent {\bf Proof.  }}

\newcommand{\Q}{\mathbb{Q}}

\newcommand{\proj}{\mathbb P}

\newcommand{\cm}{{\mathcal{M}}}

\newcommand{\al}{\alpha}

\newcommand{\be}{\beta}

\newcommand{\la}{\lambda}

\newcommand{\Aut}{\operatorname{Aut}}

\newcommand{\cmbar}{\overline{\cm}}

\newcommand{\ptz}{\partial_0}

\def\atp#1#2{\stackrel{\scriptstyle{#1}}{\scriptstyle{#2}}}

\def\wi#1#2{\langle #1 \rangle_{#2}}
\def\qed{{\hfill{\large $\Box$}}}

\begin{document}
\pagestyle{plain}
\title{The Gromov-Witten potential of a point, Hurwitz numbers,
and Hodge integrals\footnote{1991 Mathematics Subject Classification:
Primary 14H10, 81T40;
Secondary 05C30, 58D29}}
\keywords{Hurwitz numbers, Gromov-Witten potential, moduli space,
ramified covers, recurrences, Itzykson-Zuber ansatz, combinatorialization,
Hodge integrals}

\author{I.P. Goulden}
\address{Dept. of Combinatorics and Optimization,
University of Waterloo, Waterloo, Ontario, Canada}
\thanks{Supported by a research grant from the Natural Sciences
and Engineering Research Council of Canada}
\email{ipgoulden@math.uwaterloo.ca}

\author[D.M.Jackson]{D.M. Jackson}
\address{Dept. of Combinatorics and Optimization,
University of Waterloo, Waterloo, Ontario, Canada}
\thanks{Supported by a research grant from the Natural Sciences
and Engineering Research Council of Canada}
\email{dmjackson@math.uwaterloo.ca}

\author{R. Vakil}
\address{Dept. of Mathematics, MIT, Cambridge, MA~02139}
\email{vakil@math.mit.edu}
\thanks{Partially supported by a research grant DMS~9970101 from the
National Science Foundation}

\date{September 30, 1999.}

\begin{abstract}
Hurwitz numbers, which count certain covers of the projective line
(or, equivalently, factorizations of permutations into transpositions),
have been extensively studied for over a century.  The
Gromov-Witten potential $F$ of a point, the generating series for
Hodge integrals on the moduli space of curves, has been a central
object of study in Gromov-Witten theory.  We define a slightly
enriched Gromov-Witten potential $G$ (including integrals involving
one ``$\la$-class''), and show that after a non-trivial change of variables,
$G=H$ in positive genus, where $H$ is a generating series
for Hurwitz numbers.  We
prove a conjecture of Goulden and Jackson on higher genus Hurwitz
numbers, which turns out to be an analogue of a genus expansion ansatz
of Itzykson and Zuber.   As consequences, we have new combinatorial
constraints on $F$, and a much more direct proof of the ansatz of Itzykson
and Zuber.

We can produce recursions and explicit formulas for Hurwitz numbers;
the algorithm presented should prove ``all'' such recursions.
Furthermore, there are many more recursions than previously suspected
from geometry (and indeed they should exist in all genera); as
examples we present surprisingly simple new recursions in genus up to 3
that are geometrically mysterious.

As we expect this paper also to be of interest to combinatorialists, we have
tried to make it as self-contained as possible, including reviewing some
results and definitions well known in algebraic and symplectic geometry,
and mathematical physics.
\end{abstract}

\maketitle
\tableofcontents

%--------------------------------------------------------
{\parskip=12pt	% the closing bracket for this is at the end of the
		% body of text, just before the bibliography.

\section{Introduction}
The moduli space $\cmbar_{g,n}$ of $n$-pointed genus $g$ curves,
with stability condition
\begin{eqnarray}\label{Mstabcond}\remind{Mstabcond}
2g-2+n > 0
\end{eqnarray}
has dimension
\begin{eqnarray}\label{Mdim}\remind{Mdim}
3g-3+n. 
\end{eqnarray}
It is the Deligne-Mumford compactification of the moduli space $\cm_{g,n}$ of
smooth $n$-pointed genus $g$ curves.
It has $n$ natural line
bundles $\mathbb{L}_i$ (roughly, the cotangent space to the $i$th marked
point) and a natural rank $g$ vector bundle $\mathbb{E}$ (the Hodge bundle;
its fibers corresponds to global differentials on the curve.
% $H^0(C, \omega_C)$ over $[C, p_1, \dots, p_n]$). 
Let $\psi_i = c_1(\mathbb{L}_i)$ and $\la_k = c_k(\mathbb{E}),$
where $c_j$ is the $j$-th Chern class;
intersections of $\psi$-classes are called  {\em descendant integrals}, and
intersections of $\psi$-classes and $\la$-classes are called
{\em Hodge integrals} (see~\cite{fabp} for fuller information).

The Gromov-Witten potential $F$ of a point (Witten's total free energy
of two-dimensional gravity) is a generating series for all
descendant integrals.  Witten's conjecture (Kontsevich's theorem,
\cite{k}) and the Virasoro conjecture for a point can be expressed as the fact that $e^F$ is annihilated by
certain differential operators (see \cite{g1} for example).  We define $G$ as a 
generalization of $F$ (Section~\ref{background}), a generating
series for all intersections of $\psi$-classes and (up to) one
``$\la$-class''.  (This is part of the very large phase space of
\cite{mz}.)  Then $F$ can be easily recovered from~$G$.

Hurwitz numbers 
enumerate covers of the projective line by smooth connected curves of specified
degree and genus, with specified branching above one point, simple
branching over other specified points, and no other branching.
Equivalently, they are purely combinatorial objects counting
factorizations of permutations into transpositions
that generate a group which acts transitively on the sheets.
Hurwitz numbers have long been of interest (see, for example,
\cite{h}, \cite{v} for more recent references, and~\cite{ct} for
relation to mathematical physics).
Let $H$ be a generating series for Hurwitz numbers (defined
precisely in Section~\ref{background}).

It is straightforward (if tedious) to produce expressions for Hurwitz
numbers for any given degree (see~\cite{h} and~\cite{eehs} for degrees
up to 6), but geometrical arguments are required for obtaining
expressions for fixed genus and it is the latter that we consider.

\subsection{Recursions and Gromov-Witten theory}

One proof of the power of the theory of stable maps is the large
number of striking recursions it has produced for solutions to
classical problems in enumerative geometry, often as consequences of
``topological recursion relations''. The
original example was Kontsevich and Manin's remarkable recursion for
rational plane curves (\cite{km} Claim 5.2.1).  Eguchi, Hori, and
Xiong~\cite{ehx} used the Virasoro conjecture to find a recursion for
genus 1 plane curves (proved in \cite{p} and \cite{dz}).  Similar recursive
structure also underlies characteristic numbers in low genus
(\cite{ek},
\cite{v2}, \cite{gp}).

There are strong analogies between plane curves and covers of the
projective line.  Similar techniques in Gromov-Witten theory have
produced recursions for Hurwitz numbers (see \cite{fanp} pp. 17--18 or
\cite{v2} Section 5.11 for a summary), including a genus 2 relation
conjectured by Graber and Pandharipande and proved in \cite{gj2}.
Ionel has produced recursions using topological recursion relations
and the Virasoro conjecture (\cite{i}).  Geometers have thought that
recursions among Hurwitz numbers should be rare, and should not occur in high
genus.  Philosophically, Section \ref{ca} shows that in fact
recursions are ``thick on the ground'', and that there is an algorithm
for producing (and verifying) them.  It is expected that only a few
will have straightforward (and enlightening) geometric
explanations.  (It would be
interesting to reverse the Gromov-Witten approach and, for example, to
produce relations in the cohomology of $\cmbar_{g,n}$ using
recursions, but this does not seem to be tractable.)

Recurrences can be obtained in the more general setting of ramified
coverings of surfaces of higher genera. These were considered by
Hurwitz (\cite{h}).  When his approach is carried out by means of a
cut-and-join analysis, the resulting partial differential equation
(e.g. see Section \ref{cons1}) is, of course, identical to the one for
the sphere, although the initial conditions are different. It is then
a straightforward matter to write down the recurrence for arbitrary
ramification over infinity.
\cite{lzz} have obtained such a recurrence by other methods, although boundary
conditions were not included
(see also~\cite{lzz}~Thm.~B and~\cite{gjv}~Lemma~3.1).

%-----
\subsection{Organization of the paper}
We first show that, after a non-trivial change of variables (denoted
by $\Xi$), $G=H$ in positive genus (Theorem~\ref{change}).  Hence the
Gromov-Witten potential of a point is a purely combinatorial object in
a new way.  The proof uses a remarkable formula of Ekedahl, Lando,
Shapiro, and Vainshtein (\cite{elsv} Theorem~1.1) expressing Hurwitz
numbers in terms of Hodge integrals.  In some sense this addresses an
obstacle to dealing with descendant integrals, the fact that they ``do
not admit so easily of an enumerative application'' (\cite{g1} p. 1).
(Of course, Kontsevich's original formula (\cite{k}~p.~10) is also
combinatorial, and much more useful.)  However, the awkwardness of the
change of variables makes it difficult to transpose results between
``the world of $H$'' (involving Hurwitz numbers) and ``the world of
$G$'' (involving the moduli space of curves).

Second, we prove a generalization (Theorem~\ref{Gdil}) of an ansatz of
Itzykson and Zuber (\cite{iz}~(5.32), hereinafter the ``\cite{iz}
genus expansion ansatz'').  The philosophy behind the \cite{iz} genus expansion ansatz is 
that, for a fixed genus, starting from a finite number of descendant integrals
(involving those monomials in the $\psi$'s where each $\psi$-class appears with
multiplicity at least two), one can calculate any descendant integral
using only the string equation and the dilaton equation. 
The~\cite{iz} genus expansion ansatz algebraically encodes this fact.

Thirdly, we use this to prove a conjecture of Goulden and Jackson on
Hurwitz numbers (Theorem~\ref{gjpf}, \cite{gj2} Conjecture~1.2),
revealing it as a ``genus expansion ansatz for Hurwitz numbers''.  The
erstwhile mysterious combinatorial constants in the conjecture are
actually single Hodge integrals.

As an application, we observe that there are trivial combinatorial
recurrences on $H$, which lead to new conditions satisfied by $G$ (and
hence $F$).  It would be desirable to give a new proof of Witten's
conjecture using the combinatorics of covers of the projective line.
(Not surprisingly, this appears to be very difficult, and the authors
have made little progress in this direction.)  As a second
application, Theorem~\ref{gjpf} provides an algorithm for proving and
producing recursions for Hurwitz numbers.  We produce simple (and
surprising) new recursions in genus up to 3 as examples of the
algorithm's effectiveness.  Theorem ~\ref{gjpf} also yields
explicit formulas for Hurwitz numbers of any given genus; we
give an example~(\ref{genus3}) in genus 3.

%-----
\subsection{For combinatorialists}
Conjecture~1.2~\cite{gj2} came from a combinatorial approach to Hurwitz's
encoding of ramified covers, and the proof given here suggests that
further combinatorial questions of substance remain to be investigated
(for example, the combinatorialization of Hodge integrals).
Therefore, to make this paper more accessible to combinatorialists,
we specify the essential results that are taken without proof from
algebraic and differential geometry. These are
the stability condition~(\ref{Mstabcond})
and dimension condition~(\ref{Mdim})
for $\cmbar_{g,n},$
$\la_k = 0$ unless $0 \leq k \leq g,$ 
the convention $\la_0=1,$
the genus condition~(\ref{dimlapsi}) for the nonvanishing
of Hodge integrals,
the evaluation~(\ref{basevals}) of the base values
$\wi{\tau_0^3}0,$ $\wi{\tau_1}1$ and $\wi{\la_1}1,$ 
the string~(\ref{streq}) and dilaton~(\ref{dilatoneq}) equations for Hodge integrals,
the Riemann-Hurwitz formula~(\ref{rhformula}) for the genus of a ramified cover
and the result~(\ref{elsvf}) of Ekedahl, Lando, Shapiro and Vainshtein
relating Hurwitz numbers to Hodge integrals.  
References are given to sources
where the proofs of these are to be
found.
% All of our work with Hodge integrals is through the dilaton and
% string equation, and Faber's program~\cite{fab}.
All of our work with Hodge integrals is through the dilaton and
string equation which, in a real sense, removes the need to use
the primary definition~(\ref{tauMbar}) of Hodge integrals. 
We also use Faber's program~\cite{fab}. 

It is hoped that, for the most part, the remainder
of the paper can be read without recourse to algebraic or differential
geometry.

%------------------------------------------------------------

\section{Background}\label{background}
We begin with the necessary background on the generating series
$F, G$ and $H$ that are central to the subject of this paper.

\remind{background}

%-----
\subsection{Algebraic notation}
Suppose $\al$ is the composition $d=\al_1 + \dots + \al_m$ where the
$\al_i$ are {non-negative} integers.  Set $l(\al) = m$, the
{\em length} of $\al,$ and let $\#\Aut(\al)$ be the number of automorphisms
of the multiset $\{ \al_1, \al_2,
\dots, \al_m \}$ (so if $\be_j$ of the $\al_i$'s are $j$, then $\#
\Aut(\al) = \be_0! \be_1!$ \dots).
If the $\al_i$ are positive and non-decreasing, we write $\al \vdash d,$
and $\alpha$ is a {\em partition}.
If, furthermore, all $\al_i$ are at least~2, we write $\al \models d$.

Throughout, $t=(t_0,t_1,\ldots)$ and $p=(p_1,p_2,\ldots)$ where
$t_0,t_1,\ldots$ and $p_1,p_2,\ldots$ are indeterminates.
Thus, for example, $\mathbb{Q}[[t]]=\mathbb{Q}[[t_0,t_1,\ldots]]$ 
and $\mathbb{Q}[[x,p]]=\mathbb{Q}[[x,p_0,p_1,\ldots]]$.
If $Z$ is a polynomial in $t$, let $\left[ \frac {t_0^{k_0}} {k_0!} \dots
\frac {t_i^{k_i}} {k_i!} \right] Z$ be the coefficient of
$\frac {t_0^{k_0}} {k_0!} \dots
\frac {t_i^{k_i}} {k_i!}$ in $Z$.

Functional equations of the
form $v=xg(v),$ where $v\in\mathbb{Q}[[x]]$ and $g(0)\neq0,$
% type given in~(\ref{sequ})
have a unique solution $v(x)$ in $\mathbb{Q}[[x]]$ and an explicit expression for $f(v),$
where $f$ is an arbitrary series, may be obtained by Lagrange inversion
(see, for example~\cite{gj3} Section~1.2; also known as Lagrange's Implicit
Function Theorem). We will invoke Lagrange inversion a number of times, particularly when
deriving explicit expressions for certain Hurwitz
numbers.

%-----
\subsection{The Gromov-Witten and enriched Gromov-Witten potentials
$F$ and $G$ of a point}

Recall that
$\psi_i$ (resp. $\la_k$) is a codimension 1 (resp. $k$) Chow class on
$\cmbar_{g,n}$ where $1 \leq i \leq n$ (resp. $0 \leq k \leq g$;	
$\la_0 = 1$).  For non-negative integers $\theta_1,\ldots,\theta_n$ define 
\begin{eqnarray} \label{tauMbar}\remind{tauMbar}
\wi{ \tau_{\theta_1} \dots \tau_{\theta_n} \la_k} g = 
    \int_{\cmbar_{g,n}} \psi_1^{\theta_1} \dots \psi_n^{\theta_n} \la_k 
\end{eqnarray}
if
\begin{eqnarray}\label{dimlapsi}\remind{dimlapsi}
3g-3+n=\sum\theta_i+k
\end{eqnarray}
and $2g-2+n>0,$ and is $0$ otherwise.
(Condition~(\ref{dimlapsi}) arises  because non-zero intersections can
only occur when the sums of the codimension
of the classes intersected equals the dimension $3g-3+n$ of the
space $\cmbar_{g,n}.$) 
%-------------------
The condition equivalent to~(\ref{dimlapsi}) for
$\wi{\tau_0^{b_0}\tau_1^{b_1}\ldots\la_k}g$ is
\begin{eqnarray}\label{dimlapsi2}\remind{dimlapsi2}
k=\sum (1-i)b_i+3g-3.
\end{eqnarray}
%-------------------
In sums involving Hodge integrals it is convenient to include $k$
as a summation index, but then to recall
that the condition (either~(\ref{dimlapsi}) or~(\ref{dimlapsi2})) on $k$ is implicit.
When $k=0$, this agrees with the usual definition.  In particular,
\begin{eqnarray}\label{basevals}\remind{basevals}
\wi {\tau_0^3}0 = 1,\qquad \wi{\tau_1}1 = \wi{\la_1}1 = \frac 1 {24}.
\end{eqnarray}

\begin{definition} Let $g\ge0.$
The {\em genus $g$ Gromov-Witten potential} of a point is
$$ F_g(t) =\sum_{n \geq 0} \frac 1 {n!}
\sum_{\theta_1,\dots,\theta_n\geq 0}
t_{\theta_1}\dots t_{\theta_n}\wi{\tau_{\theta_1}\dots\tau_{\theta_n}} g.  $$
where the sum is constrained by~(\ref{dimlapsi}) with $k=0.$

\noindent The {\em Gromov-Witten potential of a point} is 
$$F=  \sum_{g \geq 0} y^{g-1} F_g.$$

\noindent The {\em genus $g$ enriched  Gromov-Witten potential} of a point is
\begin{equation*} \label{Gdef}
G_g(t) =  \sum_{n \geq 0} \frac 1 {n!}
\sum_{\theta_1,\dots,\theta_n\geq 0, 0 \leq k\leq g }(-1)^k 
t_{\theta_1} \dots t_{\theta_n} \wi{ \tau_{\theta_1} \dots \tau_{\theta_n}\la_k} g.
\end{equation*} \remind{Gdef}
where the sum is constrained by~(\ref{dimlapsi}).

\noindent The {\em enriched Gromov-Witten potential} of a point is
$$G = \sum_{g\ge0} G_g y^{g-1}.$$
\end{definition}

It will be convenient to use $G_g$ in the form
\begin{eqnarray*}
G_g(t)=\sum_{a_1,a_2,\ldots\ge0, 0\le k\le g} (-1)^k
\wi{\tau_0^{a_0}\tau_1^{a_1}\ldots\lambda_k}g
\frac{t_0^{a_0}}{a_0!}\frac{t_1^{a_1}}{a_1!}\ldots
\end{eqnarray*} 
where the sum is constrained by~(\ref{dimlapsi2}).
(The $(-1)^k$ in the definition of $G_g$ is included to
make the change of variables simpler.)
Note that $F_0=G_0.$
Note also that $F$ can be recovered from $G$ by substituting $v^{1-i} t_i$ for $t_i$,
and $v^3 y$ for $y$, and letting $G^\#(t,y,v)$ be the
resulting generating series in the $t_i$, $y$, and $v$. 
Then $F(t,y) = G^\#(t,y,0)$ and $G(t,y) = G^\#(t,y,1)$. 
Phrased differently, if $t_i$ is given degree $1-i$ and $y$ is given degree
$3,$  then $G_g$ has terms only in degrees $0$ to $g,$ and $F_g$ is the
degree $0$ part of $G_g.$
Also, 
$$
\left[ \frac {t_0^{l_0}} {l_0!} \dots
\frac {t_i^{l_i}} {l_i!}  v^k \right] G^\#_g
= (-1)^k \wi{ \tau_0^{l_0} \dots \tau_i^{l_i} \la_k }g.$$

The following equations facilitate the systematic elimination of $\tau_0$ and
$\tau_1$ from the Hodge integrals.
Let $a_0,a_1,\ldots$ be non-negative integers.
The {\em string equation} (or puncture equation) is
\begin{eqnarray}\label{streq}\remind{streq}
\wi{\tau_0^{a_0+1}\tau_1^{a_1}\dots\la_k}g= \sum_{i \geq 0} a_{i+1}
\wi{\tau_0^{a_0}\dots\tau_i^{a_i+1}\tau_{i+1}^{a_{i+1}-1}\dots\la_k}g, 
\end{eqnarray}
unless $g=0$, $k=0$, $a_0=2$, and all other $a_i$ are zero (in which case the
left hand side is $\wi{\tau_0^3}0=1$ by~(\ref{basevals})).
In genus $0,$ for example,
\begin{eqnarray}\label{eMbar0n}
\int_{\cmbar_{0,n}} \psi_1^{\theta_1} \dots \psi_n^{\theta_n}= \binom {n-3}
{\theta_1,\ldots,\theta_n}
\end{eqnarray}
by a trivial induction from the string equation
(observe that one of the $\theta_i$ has to be zero, so the string equation may be
applied) with $\wi{\tau_0^3}0=1$ as the base case.

The {\em dilaton equation} is
\begin{equation}\label{dilatoneq}\remind{dilatoneq}
\wi{\tau_0^{a_0}\tau_1^{a_1+1}\tau_2^{a_2}\ldots\la_k}g=\left(2g-2+\sum_i a_i\right)
\wi{\tau_0^{a_0}\tau_1^{a_1}\tau_2^{a_2}\ldots\la_k}g,
\end{equation}
unless $g=1$, $k=0$, and $a_i$ are all zero
(in which case the left hand side is $\wi{\tau_1}1=1/24$ by~(\ref{basevals})). 
The proofs of the string and dilaton equations are the same as the usual
proofs~(for example, \cite{l}~p.~191) when
no $\la$-class is present so we suppress them.  In particular, by induction, we
obtain the following repeated form of the dilaton equation from the
dilaton equation:
if $a=a_0+a_1 + \cdots$, then
\begin{equation}\label{repeateddilaton} \remind{repeateddilaton}
\wi {\tau_0^{a_0} \tau_1^{a_1} \tau_2^{a_2} \dots \la_k }g
= \frac { (a+2g-3)!} { (a+2g-3-a_1)!} \wi{\tau_0^{a_0} \tau_2^{a_2} \dots \la_k}g
\end{equation}
(except when the equation does not make sense, i.e.  when
$g=0$ and $a-a_1<3$, or $g=1$ and $a-a_1 = k=0$),
expressing the consequence of eliminating each $\tau_1.$
The string and dilaton equations can be easily translated into
differential equations for $G_g.$

%-----
\subsection{The Hurwitz generating series $H$}

Fix a genus $g$, a degree $d$, and a partition
$(\al_1,\ldots,\al_m)$ of $d$ with $m$ parts.  Let
\begin{eqnarray}\label{rhformula}\remind{rhformula}
r = d+m+2(g-1),
\end{eqnarray}
so a branched cover of $\proj^1$, with monodromy above
$\infty$ given by $\al$, and $r$ other specified simple branch points
(and no other branching) has genus $g$ (by the Riemann-Hurwitz formula). 
Let $H^g_{\al}$ be the number of such branched covers that are connected. 
(We do not take the branched points over $\infty$ to be labelled.)  

Ekedahl, Lando, Shapiro and Vainshtein have announced a remarkable
formula~(\cite{elsv}~Theorem~1.1)
\begin{equation}\label{elsvf}\remind{elsvf}
H^g_\al = \frac {r!} { \# \Aut(\al)}
 \prod_{i=1}^m \frac {{\al_i}^{\al_i}} {\al_i!} 
\int_{\cmbar_{g,m}}      \frac { 1-\la_1 + \dots \pm \la_g} {\prod (1-\al_i \psi_i)}
\end{equation}
that expresses Hurwitz numbers in terms of Hodge integrals.
\remind{elsvf} (In their statement, $\Lambda_{g;n}$ should be replaced
by $\Lambda_{g;n}^\vee$; $\Lambda_{g;n}$ is the Hodge bundle $\mathbb{E}.$) 
One should be cautious before using this
formula when $(g,m)=(0,1)$ or $(0,2)$ (as the moduli functor
$\cmbar_{g,m}$ is not a Deligne-Mumford stack), but we will
not use these two degenerate cases.

A proof of ~(\ref{elsvf}) using virtual localization~(\cite{gp1}) in
the moduli space of stable maps to $\proj^1$ will appear in \cite{gv}.
It will also be explained there how (\ref{elsvf}) would follow
quickly from virtual localization on an appropriate ``relative''
moduli space, not yet defined in the algebraic category (yielding
relative Gromov-Witten invariants; see
\cite{lr} Section 7 and \cite{ip} for discussion in the symplectic
category, and \cite{gathmann} for some discussion in the algebraic
category in the case $g=0$).  In the case where there is no
ramification above $\infty$ (i.e.  $\al = (1^d)$), the argument
reduces to Fantechi and Pandharipande's independent proof of
(\ref{elsvf}), \cite{fanp}~Theorem~2.

\begin{definition} 
The {\em Hurwitz generating series} is
$$H = \sum_{g \geq 0} H_g y^{g-1},$$
where $H_g$ is the generating series
$$H_g=H_g(x,p)=\sum_{d\ge1,\al\vdash d}\frac{H^g_{\al}}{r!}p_{\al}x^d$$
for the $H^g_\al,$ $p_{\al}$ and $x$ are indeterminates, and where
$2-2g = d-r+l(\al).$
\end{definition}
Note that $e^H$ counts all covers, not just connected ones.
($H_g$ is denoted by $F_g$ in~\cite{gj2}.)

Goulden and Jackson have conjectured that $H_g$ is of
a particular form in terms of an implicitly defined set of variables
$\{\phi_i(s,p)\colon i\ge0\}$ defined as follows. Let
\begin{equation} \label{defphi}\remind{defphi}
\phi_i(z,p) = \sum_{n \geq 1} \frac {n^{n+i}} { n!} p_n z^n,
\end{equation}
where $i$ is an integer,
be a formal power series (called $\psi_i(z,p)$ in \cite{gj2}).
Then, through the functional equation
\begin{equation} \label{sequ}\remind{sequ}
s = x e^{\phi_0(s,p)},
\end{equation}
$s$ is uniquely defined as a formal power series in $x$ (and $p$).  

In particular, $H_0$ and $H_1$ are given in~(\ref{H0xp}) and~(\ref{H1xp}),
respectively. The remaining $H_g$ are the subject of the following conjecture.

\begin{conjecture} [Goulden-Jackson \cite{gj2}~Conj. 1.2] \label{gjconj}
For $g \geq 2$,
\begin{equation} \label{gjconjf}
H_g(x,p) = \sum_{e = 2g-1}^{5g-5} \frac 1 { (1-\phi_1(s,p))^e }
\sum_{n=e-1}^{e+g-1}
\!\!\!\!\!
\sum_{  \substack { {\theta \models n} \\
{l(\theta) = e-2(g-1)}}}
\!\!\!\!\!
\frac{K_{\theta}^g}{\#\Aut(\theta)} 
{\phi_{\theta_1}(s,p)}
{\phi_{\theta_2}(s,p)}
\dots
\end{equation}
for some rational numbers $K^g_\theta$.
\end{conjecture}
\remind{gjconjf}
\remind{gjconj}

We prove this conjecture (Theorem~\ref{gjpf}).  Remarkably, each
unknown constant $K^g_\theta$ turns out to be a single Hodge integral, up to sign.  

\begin{remark}\label{gjdisc} \remind{gjdisc}
{\rm 
Goulden and Jackson proved Conjecture \ref{gjconj} for $g=2$, and conjectured
explicit values for certain $K^g_\theta$ (for $g=3$ and all $\theta$
\cite{gj2} Appendix A, and for $(e,l(\theta)) = (2g-1, 1)$
and all admissible $g$ and $n$ \cite{gj2} p. 3); we discuss these further
in Section~\ref{gjpfcor}.  
}
\end{remark}

%-----
\subsection{The relationship between $H_g$ and $G_g$}
The following is a useful result that connects  $H_g$ and $G_g.$
Throughout this section and the next we will make use of the mapping
$$\Xi\colon t_k\longmapsto\phi_k(x,p),$$
for $k\ge0,$
extended as a homomorphism to $\mathbb{Q}[[t]].$
\begin{theorem}\label{change}
If $g>0$, then $H_g(x,p) = \Xi\,G_g(t).$
\end{theorem}
\remind{change}
\bpf
For $g>0$, by (\ref{elsvf}),
\begin{eqnarray*}
H_g &=& \sum_{\al \vdash d}  \frac 1 {\# \Aut(\al)} \frac {\prod 
\al_i^{\al_i}} {\prod \al_i!} p_{\al} x^d \int_{\cmbar_{g,m}}
\frac {1 - \la_1  +  \dots \pm \la_g} { \prod (1-\al_i \psi_i)}  \\
&=& \sum_{\al_1 + \dots + \al_m = d} \frac 1 {m!} \frac {\prod 
\al_i^{\al_i}} {\prod \al_i!} p_{\al} x^d \int_{\cmbar_{g,m}}
\frac {1 - \la_1  +  \dots \pm \la_g} { \prod (1-\al_i \psi_i)} \\
 &=& \sum_m \frac 1 {m!} \sum_{\al_1, \dots, \al_m \geq 1} 
\prod_i \left( \frac { \al_i^{\al_i} p_{\al_i} x^{\al_i}} {\al_i!} \right) \\
&\mbox{}& \cdot \sum_{ \atp { b_1 + \dots + b_m = 3g-3+m-k}  {0 \leq k \leq g, b_i \geq 0} } 
\int_{\cmbar_{g,m}} (\al_1 \psi_1)^{b_1} \dots (\al_m 
\psi_m)^{b_m} (-1)^k \la_k \\
&=& \sum_m \frac 1 {m!} \sum_{ \atp { b_1 + \dots + b_m = 3g-3+m-k}  {0 \leq k \leq g, b_i \geq 0} } 
(-1)^k \wi{ \tau_{b_1} \dots \tau_{b_m} \la_k } g    \\
&\mbox{}& 
\cdot\sum_{\al_1, \dots, \al_m \geq 1}  
\prod_i \left( \frac { \al_i^{\al_i+b_i} p_{\al_i} x^{\al_i}} {\al_i!} \right).
\end{eqnarray*}
Hence
$$H_g = \sum_{m \geq 0} \frac 1 {m!} 
\sum_{b_1, \dots, b_m \geq 0, 0\le k\le g} (-1)^k 
\left( \prod_{i=1}^m \phi_{b_i} (x,p) \right) \wi{\tau_{b_1} \dots \tau_{b_m} 
\la_k} g .  $$
The result then follows from (\ref{Gdef}). \qed

If $g=0,$ the above statement must be modified. The formula~(\ref{elsvf})
applies when $l(\al)\ge3,$ so if $H_g[m]$ is the summand of $H_g$ 
corresponding to all $\al$ with $l(\al)=m,$ then
$$H_0 = H_0[1]+H_0[2]+\sum_{m\ge3}H_0[m] = H_0[1]+H_0[2]+\Xi\,G_0,$$
so
\begin{eqnarray*}\label{H0expression}\remind{H0expression}
H_0 = H_0[1]+H_0[2]+\Xi\,F_0.
\end{eqnarray*}

de Jong has pointed out that the change of variables $\Xi$ is not
invertible.  In other words, ignoring the irrelevant variable $x$
by setting it equal to $1,$ $\Xi$ is not invertible.
To see this, let $\rho\colon p_n\longmapsto np_n$ and
$\sigma\colon t_n\longmapsto t_{n+1}.$ Then $\rho\Xi=\Xi\sigma.$
But $\rho$ is invertible and $\sigma$ is not. Thus $\Xi$ is
not invertible.

%--------------------------------------------------------------------

\section{Structure theorems for $G$ and $H$}

For $k\ge0,$ let
\begin{eqnarray}\label{defIk}\remind{defIk}
I_k = \sum_{i\ge0} t_{k+i} \frac{I_0^i}{i!}.
\end{eqnarray}
When $k=0,$ this is a functional equation that, by Lagrange inversion,
uniquely defines $I_0\in\mathbb{Q}[[t]],$ and thence $I_k$ is uniquely defined as
a series in $\mathbb{Q}[[t]]$ for all $k\ge0.$
If $t_0=0,$ the unique solution of~(\ref{defIk}) is $I_0=0,$ so that with this specialization 
\begin{eqnarray}\label{icondIk}\remind{icondIk}
I_k=t_k \mbox{ for } k\ge1.
\end{eqnarray}

%------
\subsection{Structure theorem for $G$}
The following is a generalization of the \cite{iz} genus expansion ansatz.
This argument also gives a much more direct proof of
the original \cite{iz} genus expansion ansatz, by
``setting $\la_k=0$'' for $k>0$ (excising terms for all $\theta$ such that
$\sum_j(1-j)\theta_j+3g-3>0$).
(The only proof of the \cite{iz} ansatz in the literature
known to the authors is in \cite{eyy}.)
 Denote $\partial/\partial t_i$ by $\partial_i$
for the sake of brevity.

\begin{theorem} [Genus expansion ansatz]\label{Gdil}
If $g>1$, 
\begin{eqnarray} \label{Gdileq} 
G_g(t) &=& \frac1{(1-I_1)^{2g-2}}G_g\left(0,0,\frac{I_2}{1-I_1},\frac{I_3}{1-I_1},\ldots\right) \\
&=& \sum_{      \substack{      {\sum_{2 \leq j \leq 3g-2} (j-1) l_j} \\    { + k= 3g-3 }  }     }
(-1)^k \frac{\wi {\tau_2^{l_2}\tau_3^{l_3} \dots \tau_{3g-2}^{l_{3g-2}} \la_k }g}
 { (1-I_1)^{2(g-1)+\sum l_j}}
{ \frac{I_2^{l_2}}{l_2!}  \dots \frac{I_{3g-2}^{l_{3g-2}}} {l_{3g-2}!} }.
\label{Gdileq2}
\end{eqnarray}
\remind{Gdileq}
\end{theorem}
\remind{Gdil} 
% local def
\newcommand{\Gg}{G_g(0,I_1,I_2,\ldots)}
(It is straightforward to show that the right sides of equations
(\ref{Gdileq}) and (\ref{Gdileq2}) are the same.) 

In \cite{fabp2} Section 2.1, Faber and Pandharipande use the
terminology ``primitive'' to denote Hodge integrals without $\tau_0$
or $\tau_1$.  Essentially the formal derivation here (like the work of
\cite{iz}) is to write an explicit formula for $G_g$ in terms of
primitive Hodge integrals.  Viewed in this way, it is clear there are
only finitely many degrees of freedom for each genus (as there are
only finitely many primitive Hodge integrals for a fixed genus); the
interesting part is the precise form.
 
\bpf  Let $\Delta = \sum_{m \geq 0} t_{m+1} \partial_m - \ptz.$
Then, from the string equation~(\ref{streq}), $$\Delta G_g(t)=0,$$  for $g>0,$ 
and $G_g(t)$ is the unique such series with the initial value
$G_g(0,t_1,\ldots)$ at $t_0=0.$ We begin the proof by exploiting this uniqueness
to establish that
\begin{eqnarray}\label{GgtGgI}\remind{GgtGgI}
G_g(t)=\Gg, \mbox{ for } g>0.
\end{eqnarray}

Let $\zeta_i=0$ if $i<0$ and $1$ if $i\ge0.$  Then, from~(\ref{defIk}), for $m,k\ge0,$
\begin{eqnarray*}
\partial_m I_k = \zeta_{m-k} \frac{I_0^{m-k}}{(m-k)!}+
\left(\sum_{i\ge1} t_{k+i} \frac{I_0^{i-1}}{(i-1)!}\right)\,\partial_m I_0,
\end{eqnarray*}
so
\begin{eqnarray*} 
\partial_m I_k = \zeta_{m-k} \frac{I_0^{m-k}}{(m-k)!}+I_{k+1}\,\partial_m I_0.
\end{eqnarray*}
Then, substituting  $k=0$ above, we obtain for $m\ge0$
\begin{eqnarray*}
\partial_m I_0 = \frac1{m!}\frac{I_0^m}{1-I_1},
\end{eqnarray*}
so, for $k,m\ge0,$
\begin{eqnarray}\label{dmIk2}\remind{dmIk2}
\partial_m I_k = \zeta_{m-k} \frac{I_0^{m-k}}{(m-k)!}+ \frac{I_0^m}{m!}\frac{I_{k+1}}{1-I_1}.
\end{eqnarray}
Now, by the chain rule,
\begin{eqnarray*}
\Delta\Gg 
&=& \sum_{k\ge1} \left(\sum_{m\ge0} t_{m+1} \partial_m I_k -\ptz I_k \right) \frac{\partial}{\partial I_k}\Gg.
\end{eqnarray*}

But, from~(\ref{dmIk2}),
\begin{eqnarray*}
\sum_{m\ge0}t_{m+1}\partial_m I_k -\ptz I_k 
&=&\sum_{m\ge k} t_{m+1}\frac{I_0^{m-k}}{(m-k)!} + 
\frac{I_{k+1}}{1-I_1} \sum_{m\ge0} t_{m+1}\frac{I_0^m}{m!} - \frac{I_{k+1}}{1-I_1} \\
&=& 0,
\end{eqnarray*}
for $k\ge1.$ Thus $\Delta \Gg=0.$
But $\Gg \vert_{t_0=0} = G_g(0,t_1,t_2,\ldots)$
from~(\ref{icondIk}), and thus we have established~(\ref{GgtGgI}) by the uniqueness argument.

To complete the proof, we use the repeated form~(\ref{repeateddilaton})
of the dilaton equation for $g>1.$ 
\begin{eqnarray*}
\Gg &=& \sum_{b_1,b_2,\ldots\ge0}
(-1)^{\sum_{i\ge1}(1-i)b_i+3g-3}
\wi{\tau_1^{b_1}\tau_2^{b_2}\ldots \lambda_k}g
\frac{I_1^{b_1}}{b_1!}\frac{I_2^{b_2}}{b_2!}\cdots \\
&=&  \sum_{b_2,b_3,\ldots\ge0}
(-1)^{\sum_{i\ge2}(1-i)b_i+3g-3}
\wi{\tau_2^{b_2}\tau_3^{b_3}\ldots \lambda_k}g
\frac{I_2^{b_2}}{b_2!}\frac{I_3^{b_3}}{b_3!}\cdots \\
&\mbox{}& \cdot\sum_{b_1\ge0}
\binom{-(b_1+b_2+\cdots)-2g+2}{b_1} I_1^{b_1}
\end{eqnarray*}
from~(\ref{repeateddilaton}). Thus
$$\Gg=\frac1{(1-I_1)^{2g-2}}G_g\left(0,0,\frac{I_2}{1-I_1},\frac{I_3}{1-I_1},\ldots\right),\mbox{ for } g>1,$$ 
and the result now follows from~(\ref{GgtGgI}). \qed

%---------------------------------------------------------------------
\subsection{Structure theorem for $H$}

%-----
% \subsection{The main theorem}
We now give the main structure theorem for $H.$
\begin{theorem} [\cite{gj2} Conjecture~1.2]  \label{gjpf}
Conjecture \ref{gjconj} is true, with
\begin{equation}
\label{specialk}
K^g_{\theta} = (-1)^k \wi {\tau_{\theta_1} \tau_{\theta_2} \dots \la_k}g,
\end{equation}
where $k=\sum_j (1-j) \theta_j + 3g-3$.
\end{theorem}
\remind{gjpf}

\bpf
{}From Theorem~\ref{change} with $g>0,$ $H_g(x,p)=\Xi\,G_g(t)$ where,
from Theorem~\ref{Gdil} (\ref{Gdileq2}), for $g>1,$

\begin{equation*}
G_g = \!\!\!\!\!\!\!\!\! 
\sum_{\sum_{2 \leq j \leq 3g-2} (j-1) l_j + k= 3g-3}
\!\!\!\!\!\!\!\!\! \!\!\!\!\!\!\!\!\! 
(-1)^k \frac{\wi {\tau_2^{l_2}\tau_3^{l_3} \dots \tau_{3g-2}^{l_{3g-2}} \la_k }g}
 { (1-I_1)^{2(g-1)+\sum l_j}}
{ \frac{I_2^{l_2}}{l_2!}  \dots \frac{I_{3g-2}^{l_{3g-2}}} {l_{3g-2}!} }.
\end{equation*}
We want to prove~(\ref{gjconjf}), for $g \geq 2;$ that is,
\begin{equation*}
H_g(x,p) = \sum_{e = 2g-1}^{5g-5} \frac 1 { (1-\phi_1(s,p))^e }
\sum_{n=e-1}^{e+g-1} \sum_{  \substack { {\theta \models n} \\
{l(\theta) = e-2(g-1)}}}
\frac{K_{\theta}^g}{\#\Aut(\theta)} 
{\phi_{\theta_1}(s,p)}
{\phi_{\theta_2}(s,p)} \dots
\end{equation*}
where $K^g_{\theta}$ satisfies (\ref{specialk}).
Since this can be rewritten in the form
$$ H_g(x,p)= \!\!\!\!\!\!\!\!\!
\sum_{\sum_{2 \leq j \leq 3g-2} (j-1) l_j + k= 3g-3}
\!\!\!\!\!\!\!\!\! % \!\!\!\!\!\!\!\!\!
\frac{K^g_{(2^{l_2}3^{l_3}\ldots)}}
{(1-\phi_1(s,p))^{2(g-1)+\sum l_j}}
\frac{\phi_2(s,p)^{l_2}}{l_2!}\cdots\frac{\phi_{3g-2}(s,p)^{l_{3g-2}}}{l_{3g-2}!},
$$
the proof is therefore complete if we can establish that
$\Xi\,I_k(t)=\phi_k(s,p)$ for $k\ge1,$ thereby making the identification
$K_\theta^g=(-1)^k \wi {\tau_{\theta_1} \tau_{\theta_2} \dots \la_k}g.$

From~(\ref{defphi}) and~(\ref{sequ}), for $k\ge0,$
\begin{eqnarray*}
\phi_k(s,p) &=& \sum_{n\ge0} \frac{n^{n+k}}{n!} p_n x^n e^{n\phi_0(s,p)} \\
&=& \sum_{m,n\ge0} \frac{n^{n+k+m}}{n!} p_n x^n \frac{\phi_0(s,p)^m}{m!},
\end{eqnarray*}
so
\begin{eqnarray*}\label{phiksp}\remind{phiksp}
\phi_k(s,p) &=& \sum_{m\ge0}\phi_{k+m}(x,p)\frac{\phi_0(s,p)^m}{m!}.
\end{eqnarray*}
By comparing this with the definition~(\ref{defIk}) of $I_k,$
it follows that $\Xi\,I_k(t)=\phi_k(s,p)$ for $k\ge0,$ completing the proof. \qed

We record the observation on the action of $\Xi$ that
\begin{eqnarray*}\label{XiIphi}\remind{XiIphi}
\Xi\,I_k=\phi_k(s,p), \mbox{ for } k\ge0.
\end{eqnarray*}
Thus we have established the connexion between the indeterminates
$x,p_i$ on the Hurwitz side and the indeterminates $t_r$ and $I_r$ on the
Gromov-Witten side (see Section \ref{door}).

%-----
\subsection{Analogous statements in genus 0 and 1}
We note that~(\cite{gj0}~Proposition~3.1(1))
\begin{eqnarray}\label{H0xp}\remind{H0xp}
\left(x\frac\partial{\partial x}\right)^2H_0(x,p)=\phi_0(s,p).
\end{eqnarray}
In the light of Theorem~\ref{change}, stating that $\Xi\,G_g(t)=H_g(x,p)$ for $g>0$,
earlier statements in geometry and in combinatorics can now be seen to be
equivalent.
In genus 1,
\begin{eqnarray}\label{H1xp}\remind{H1xp}
H_1(x,p) = \Xi\,G_1(t) = \frac 1 {24} \left( \log (1 - \phi_1(s,p))^{-1}
- \phi_0(s,p) \right)
\end{eqnarray}
(\cite{v}, \cite{gj1} Theorem~4.2), 
and
$$\Xi\,F_1(t) = \frac 1 {24} \log (1 - \phi_1(s,p))^{-1}$$
(\cite{iz}~(5.30), \cite{eyy}~(3.7), \cite{dw}). 
The difference $-\frac 1 {24} \phi_0(s,p)$ can be seen to be
the contribution to $\Xi\,G_1(t)$ from $\la_1$.

Surprisingly, the picture is least clear in genus 0.
$F_0(t)=G_0(t)$, and
the difference $H_0(x,p)-\Xi\,G_0(t)$ arises from where (\ref{elsvf}) breaks down:
it is a generating series for covers of $\proj^1$ with at most~2
pre-images of $\infty,$ $H_0[1](x,p)+H_0[2](x,p).$
By~\cite{gj0} or~\cite{denes},
$$H_0[1](x,p) = \phi_{-2}(x,p).$$  By \cite{a} or \cite{gj0},
$$H_0[2](x,p) = \sum_{i, j \geq 1}
\frac { (i+j-1)!} { (i-1)! (j-1)!} i^{i-1} j^{j-1} p_i p_j x^{i+j}.$$
{}From~(\ref{H0expression}), $\Xi\,F_0(t)+H_0[1](x,p)+H_0[2](x,p)=H_0(x,p)$
so, using formula (\ref{eMbar0n}) for $F_0$ and \cite{gj0}
Theorem~1.1 for $H_0,$ this gives an explicit relation.
However, it does not seem enlightening.

%-----
\begin{remark}\label{gjpfcor} \remind{gjpfcor}  
{\rm
Using Theorem~\ref{gjpf}, it follows that the conjectures of Goulden
and Jackson described in Remark~\ref{gjdisc} are true.  The
conjectured values of $K^3_\theta$ can be checked using Faber's
program~\cite{fab}.   The conjectured values of $K^g_\theta$ for $e=2g-1$,
$l(\theta)=1$ (involving coefficients of $\left( \frac {z/2} { \sin(z/2)} \right)^{k+1}$)
turn out to be equivalent to \cite{fabp} Theorem~2 and
\cite{elsv} Theorem~1.2.
}
\end{remark}

%------------------------------------------------------------------------
\section{Consequences and applications}
\label{ca}

\subsection{Combinatorial comments on Hodge integrals}
The terms that appear in Conjecture~\ref{gjconj} can be
given, in principle, a combinatorial interpretation.
The left hand side already has a combinatorial interpretation,
through Hurwitz's encoding, in terms of transitive ordered factorizations
of permutations into transpositions.

For the right hand side, $n^{n+i}$ is the number of rooted (vertex-)
labelled trees  with $i+1$ marked vertices (vertices may be multiply
marked). The generating series for this number is $\phi_i(z,p),$ where
$p_n$ records the number of vertices in a tree.
$\phi_0(z,p)$ is therefore the number of rooted labelled trees with
exactly one marked vertex. Similar interpretations can therefore be
given to $s$ and $1/(1-\phi_1(s,p))^e.$ The right hand side therefore
has an interpretation in terms of structures obtained by gluing together
and ordering collections of rooted labelled trees with marked vertices.
This suggests that $K^g_\theta,$ which has been identified up to sign as
a Hodge integral through Theorem~\ref{gjpf}, can be defined purely
combinatorially, provided the
mapping between the structures corresponding to the left hand and
right hand sides of (\ref{gjconjf}) is made explicit.  In particular,
this would involve determining how markers attached to the vertices of
the trees from the right hand side encode transitive ordered factorizations of
permutations into transpositions, that occur on the left hand side of
(\ref{gjconjf}). 
This is, of course, where the difficulty lies since the theorem itself
provides no information about the elementwise action of such a mapping.
% We therefore conclude that
% this programme could, in principle, lead to a purely algebraic proof
% of Witten's conjecture. 

%-----
\subsection{Consequences of Theorem~\ref{change}}\label{cons1} \remind{cons1}
Theorem~\ref{change} gives a new combinatorial structure on $G$
(and hence $F$), and one could hope to prove results about $F$ using
$H$, i.e. the combinatorics of branched covers.  For example, there is
a simple differential operator $T$ (the ``cut-and-join'' operator)
% mentioned in Remark~\ref{gjdisc}
annihilating $e^H$, corresponding to
the interpretation of $H$ as counting factorizations of permutations
(\cite{gj0} Lemma~2.2, and independently \cite{vthesis} p.~8), defined
as follows.

Define $H^\# = H^\#(x,y,u,p)$ by substituting $x u^2$ for $x$, $y
u^2$ for $y$, and $p_i u^{1-i}$ for $p_i$ in $H_g.$
Then $H_g^\# = \sum_{d\ge0, \al \vdash d} \frac { H^g_{\al}} {r!} p_\al x^d u^r$
where $r=l(\alpha)+d+2g-2$ is the
number of simple branch points (now marked by $u$).  Let $$T =
\frac 1 2 \sum_{a,b \geq 1} \left[ (a+b) p_a p_b \frac
{\partial}{\partial p_{a+b}} + \frac 1 y ab p_{a+b}
\frac {\partial}{\partial p_{a}} 
\frac {\partial}{\partial p_{b}} 
\right]
- \frac {\partial} {\partial u}.  $$ Then $T e^{H^\#} = 0$, and $H^\#$
is uniquely determined by this equation and the condition
$H^\#(x,y,0,p) = p_1 x$ (i.e. there is only one cover of $\proj^1$
unbranched away from $\infty$).

Note that, even the string equation becomes mysterious when translated to a
statement about $H$: $$
\frac {\partial} {\partial t_0} H = \frac 1 2 t_0^2  + x \frac \partial
{\partial x} H.
$$
It is not combinatorially clear why this should be true. 

% Still, some insight has been gained.  For example, consequences of the
% Virasoro conjecture involving Hurwitz numbers can be
% verified (Section~\ref{eleny}).

%-----
\subsection{Comments on the connexion between $H$ and $G$ (and $F$)}
\label{door}
It
is worth noting how the variables used by physicists to study $F$ 
(and that are equally useful for $G$) have exactly paralleled the
variables used by combinatorialists to study $H$.  Specifically, 
physicists (and geometers) write $F$ in terms of: 
\begin{enumerate}
\item[P1.] The variables $t_i$; $F_g, G_g \in \Q[[t]]$ are naturally generating series for all Hodge integrals.
\item[P2.] For $g>1$, $F_g$ and $G_g$ lie in a much smaller ring.  Via the genus reduction ansatz, 
Theorem \ref{Gdil}, $F_g$ and $G_g$ can be rewritten as
elements of $\Q[ 1/(1-I_1), I_2, I_3, \dots]$, and this
representation is particular simple (as only a finite number of
monomials appear, and their coefficients are each single Hodge
integrals).
\item[P3.]  It is often physically enlightening (\cite{iz}, \cite{eyy}) to rewrite the above
in terms of other variables.  Let $u_0 = \ptz^2 F_0$.  Then for $g>1$,
$$F_g, G_g \in\Q[ 1/ \ptz u_0, \ptz u_0, \ptz^2 u_0, \dots]$$ (and in fact $F_g$
has a particular bigrading
in terms of these variables, where $\deg \ptz^r u_0 = (1,r-1)$).
In \cite{eyy}, these variables are used in the proof of the \cite{iz} genus reduction ansatz.
It is not hard to translate between the $\ptz^r u_0$ and the $I_k$;
in particular, $u_0=I_0;$ see \cite{eyy} p. 284.
\end{enumerate}
Combinatorialists write $H$ in terms of:
\begin{enumerate}
\item[C1.] The variables $x$ and $p_i$; $H_g \in \Q[[x,p]]$ is a generating series
for all Hurwitz numbers.
\item[C2.] In fact, for $g > 1$,  $H_g$ lies in a much smaller 
ring:  $$H_g \in \Q[ 1/(1-\phi_1(s,p)), \phi_2(s,p), \phi_3(s,p), \dots],$$
which via $\Xi$ is the same as P2 above.
\item[C3.] Also, $H_g$ lies in $\Q[[\phi_0(x,p), \phi_1(x,p), \dots ]]$; via
$\Xi$ this is the same as P1 above.
\end{enumerate}

%-----
\subsection{Applications of Theorem~\ref{gjpf}}

Along with techniques from \cite{gj2}, Theorem~\ref{gjpf} gives a
machine for developing and proving recurrences and explicit formulas
for Hurwitz numbers, given that the necessary Hodge integrals can be
calculated by Faber's program~\cite{fab}.  As an example,
in~\cite{gj2}, a conjectured
recursion of Graber and Pandharipande was proved using the Theorem in
genus~2 (proved there). We now give further examples.

The examples are for the case in which there is no ramification
over $\infty.$ We will refer to the corresponding numbers as {\em simple
Hurwitz numbers.} They are obtained by setting $p_1=1$  and
$p_i=0$ for $i\neq1.$ Under this specialization, $\phi_i(x,p)=x$
for all $i,$ and, from~(\ref{sequ}), $s=w$ where $w$ is the unique solution of
\begin{eqnarray*}\label{wxes}\remind{wxes}
w=xe^w,
\end{eqnarray*}
and is given explicitly by
$$w=\sum_{n\ge1}n^{n-1}\frac{x^n}{n!}.$$  
Then $H_g$ becomes
$$\widetilde{H_g} = \sum_{d\ge1} \frac{H^g_{(1^d)}}{(2d+2g-2)!} x^d,$$
the generating series for simple Hurwitz numbers.

%-----
\begin{example}[A recurrence equation for genus~$3$] \label{recg3}\remind{recg3}
{\rm
{}From a geometric perspective, ``it is not likely such
simple recursive formulas [similar to Graber-Pandharipande's formula
in genus~2, and simpler recursions in genus~0 and~1 \cite{v}~Theorem~2.7
(our intercalation)]
occur in $g \geq 3$'' (\cite{fanp} p. 18).  However, using
Theorem~\ref{gjpf}, recurrences can be obtained as follows.
Let $D=x\,d/dx.$
Then
\begin{eqnarray*}
D^2\widetilde{H_0}(x)&=&w, \\
\widetilde{H_1}(x)&=&\frac1{24}\left(\log(1-w)^{-1}-w\right),\\
\widetilde{H_2}(x)&=&\frac1{5760}\left(
\frac{4w^2}{(1-w)^4}+\frac{28w^3}{(1-w)^5}\right),\\
\widetilde{H_3}(x)&=&
{\frac {1}{80640}}\,{\frac {{w}^{2}}{\left (1-w\right )^{6}}}
+{\frac {73}{ 90720}}\,{\frac {{w}^{3}}{\left (1-w\right )^{7}}}
+{\frac {37}{5184}}\,{ \frac {{w}^{4}}{\left (1-w\right )^{8}}}\\
&\mbox{}& +{\frac {89}{5184}}\,{\frac {{w}^{5} }{\left (1-w\right )^{9}}}
+{\frac {245}{20736}}\,{\frac {{w}^{6}}{\left (1-w \right )^{10}}}.
\end{eqnarray*}
These are from~\cite{gj2}, although the final two can now be obtained
from Theorem~\ref{gjpf}, with the help of Faber's program~\cite{fab}
to compute the necessary Hodge integrals.

It is convenient to set $w=1-W^{-1}$, so $D = W^2 (W-1)d/dW$.
Then $D^n \widetilde{H_g}(x)$
is a polynomial in $W$ provided $2g-2+n>0$. 
(The resemblance to the stability condition for $\cmbar_{g,n}$
is probably not coincidental; $D$ can be interpreted as marking a 
point above a fixed general point of $\proj^1$.)
For $(g,n)=(0,1),(0,2),$
$D^n \widetilde{H_g}(x)$ is a rational series in $W.$
A number of these series are given below.
\begin{eqnarray*}
D\,\widetilde{H_0}(x) &=& (1-{W}^{-2})/2 \\
\widetilde{H_1}(x) &=& {\frac {\log (W)W-W+1}{24W}} \\
D\,\widetilde{H_1}(x) &=& \left (W-1\right )^{2}/24 \\
\widetilde{H_2}(x) &=& 
\left (W-1\right )^{2}{W}^{2}\left (-6+7\,W\right )/1440 \\
\widetilde{H_3}(x) &=& 
\left (W-1\right )^{2}{W}^{4} \\
&\mbox{}&\cdot\left (720-6696\,W+
19250\,{W}^{2}-21840\,{W}^{3}+8575\,{W}^{4}\right ) /725760.
\end{eqnarray*}

%-------------------

Various relations can be found between the $D^n \widetilde{H_g}(x)$
for  $(g,n)\neq(0,0),(1,0)$ by positing a general form for them and equating coefficients
of powers of $W$ to obtain a set of linear equations for the parameters appearing
in this form.

With the form containing the twenty six terms
$\left(D^p\widetilde{H_i}\right)\left(D^q\widetilde{H_j}\right)$ for
$p+q=4,$ $i+j=3,$ and $D^p \widetilde{H_i},$ for $i=3,$ $1\le p\le 4,$
for $i=2,$ $1\le p\le 5,$ and for $i=1,$ $1\le p\le 7,$ the null space
has dimension 11.  (We choose this form for potential recursions
because this is the form of the recursions previously produced via Gromov-Witten
theory.)  Thus further conditions on the parameters may be applied,
although it is not at all clear whether there is a geometrically
natural choice to make.  One such expression, obtained by imposing
linearity, is

\begin{eqnarray*}
2880\, \widetilde{H_3} &=&
- \left({\frac {2}{49}}\, 
- {\frac {227}{294}}\, D
+ {\frac {99845}{588}}\, D^2
\right) \widetilde{H_2} \\
&\mbox{}& 
- \left({\frac {1}{490}}\, D^2
-{\frac {11}{294}}\, D^3
+{\frac {38845}{14112}}\, D^4
-{\frac { 1225}{576}}\, D^5\right)
\widetilde{H_1}.
\end{eqnarray*}
This gives the following  explicit formula for  $H^3_{(1^d)}$
linearly in terms of $H^2_{(1^d)}$ and $H^1_{(1^d)}:$
\begin{eqnarray*}
2880\, H^3_{(1^d)}
&=&-\left (24-454\,d+99845\,{d}^{2}\right )
\binom{2d+4}{2} {\frac {H^2_{(1^d)}}{294}} \\
&\mbox{}&+ {d}^{2}
\left (-288+5280\,d-388450\,{d}^{2}+ 300125\,{d}^{3}\right )
\binom{2d+4}{4} {\frac {H^1_{(1^d)}}{5880}}.
\end{eqnarray*}

Similar recursions exist for all genera, and these may be obtained in
the same way.

} \end{example}

\begin{example}[Another recurrence equation for genus 3, of ``geometric form'']
{\rm
As another example to show how common recursions are, we give a genus 3 recursion that is of 
a potentially geometrically meaningful form:
\begin{eqnarray*}
H^3_{(1^d)} &=&  f(d) \binom d 2  H^2_{(1^d)}  +  \sum_{i+j=d}  \left(
             g(i,j) \binom {2d+2} { 2i -2}  i j H^0_{(1^i)} H^3_{(1^j)} \right. \\
 & &  \left. +   h(i,j) \binom {2d+2} { 2i}  i j H^1_{(1^i)} H^2_{(1^j)}  \right).
\end{eqnarray*}
where $f(d)$, $g(i,j)$, and $h(i,j)$ are polynomials of low degree.

Any formula coming from a divisorial relation on the space of maps
would have such a form.  Even though such a divisorial relation
should not exist, a geometrically-motivated recursion might still exist
of this form; the recursion for genus 1 plane curves of \cite{ehx} has
this property, for example.  One might hope for some geometrical understanding
from such a recursion.

The terms on the right-hand side of the equation correspond to
divisors on the space of maps.  The first term corresponds to degree
$d$ genus 2 covers where two of the $d$ points mapping to the same
point of $\proj^1$ are attached; hence the multiplicity of $\binom d
2$.  The second term corresponds to maps where the cover is a genus 0
degree $i$ cover (a general such cover has $2i-2$ branch points) and a
genus 3 degree $j$ cover (a general such cover has $2j+4$ branch
points) such that two points mapping to the same point of $\proj^1$ (one on each component)
are glued together; the multiplicity $ij$ comes from the choice of the
two points, and the multiplicity $\binom {2d+2} {2i-2}$ comes from partitioning
the branch points between the two components.  The third corresponds to maps where
the cover is a genus 1
degree $i$ cover and a 
genus 2 degree $j$ cover with a point of one glued to a point of the other; the multiplicity
calculation is similar to the second term.  These divisors might appear with various multiplicities,
given by the polynomials $f$, $g$ and $h$.

Unfortunately, many such recursions can be found (by the same method
as in Example~\ref{recg3}), even if the degrees of $f$, $g$, and $h$
are required to be small.  One such is
\begin{eqnarray*}
f(d) &=& \frac{1}{1702263010}(1532127678 d-2213123851), \\
g(i,j) &=& -\frac{2}{121590215} (760192125 ij -12054428314 i\\
       &\mbox{}& -2006745110 j +1033797958), \\
h(i,j) &=& -\frac{4}{2553394515} (798201731250 ij -217500288725 i\\
       &\mbox{}& -473678414332 j -42109762821).
\end{eqnarray*}
There seems to be no reason why this recursion should admit a geometrical
explanation.
} \end{example}

%---- 

\begin{example}[A recurrence equation for genus 2]
{\rm
The method of Example~\ref{recg3} can be applied to the genus~2 case; we
suppress the details. The linear differential equation that is satisfied is
\begin{eqnarray*}
4320\widetilde{H_2}(x)= -300 D^2\widetilde{H_1}+7\left(D^5-D^4\right)\widetilde{H_0}.
\end{eqnarray*}
The corresponding linear recurrence equation is
\begin{eqnarray*}
180 H_{(1^d)}^2 = -25d^2\binom{2d+2}{2} H_{(1^d)}^1
+7d^4(d-1)\binom{2d+2}{4} H_{(1^d)}^0.
\end{eqnarray*}

For genus 2 and 3, $H^g_{(1^d)}$ has been expressed in terms of
$H^{g-1}_{(1^d)}$ and $H^{g-2}_{(1^d)}.$
A reason this is not entirely unexpected is that
$D$ preserves the parity of the degree of polynomials in $W.$
But the degree in $W$ of $D^n H^g(x)$ is $2n+5g-5,$
and the parity of this mod~2 is the parity of $g-1$ mod~2.
Polynomials of both parities are required on the
right hand side in the posited form of the differential equation to
match terms on the left hand side. This is to be expected to persist for $g\ge2.$

}
\end{example}

%-----
\begin{example}[Recurrence equations for genus 1 and 0]
{\rm
The parity argument in the previous example
suggests that, if there is a recurrence equation, it must be of degree
(at least) two for the genus 1 case, and indeed a degree two example
is known (due to Graber and Pandharipande, \cite{v2} Section 5.11 or
\cite{fanp} p. 18).
This recurrence can be rewritten as the differential equation
$$
D \widetilde{H_1} = D^3 \widetilde{H_0} /24 - D^2 \widetilde{H_0}/24 +  
\left( D^2 \widetilde{H_0} \right) \left( D \widetilde{H_1} \right)
$$
which is an immediate consequence of the observations that 
$D\widetilde{H_1}(x)=(W-1)^2/24$,
$D^2\widetilde{H_0}(x)=1-W^{-1}$ and
and $D^3\widetilde{H_0}(x)=W-1.$

An even simpler recursion exists originating from the differential equation
$$D\widetilde{H_1}=\frac1{24}\left(D^3\widetilde{H_0}\right)^2.$$
This gives
\begin{eqnarray}
\label{genus1}
H^1_{(1^d)}=\frac1d
\binom{2d}4
\sum_{i=1}^{d-1} i^3(d-i)^3
\binom{2d-4}{2i-2}
H^0_{(1^i)}H^0_{(1^{d-i})}.
\end{eqnarray}
The differential equation is  an immediate consequence of the above
expressions for $D\,\widetilde{H_1}$ and $D^3\widetilde{H_0}.$
Although it might not be difficult to prove (\ref{genus1}) 
geometrically, there was no geometrical reason to suspect its existence.

The sphere is included for completeness from this point of view.
Again, by the parity argument, a recurrence of degree two is expected.
The simplest such differential equation is
$$D^2\widetilde{H_0} =\frac12 \left(D^2\widetilde{H_0}\right)^2
+D\widetilde{H_0},$$
which is an immediate
consequence of the observations that
$D^2\widetilde{H_0}(x)=1-W^{-1}$ and $D\widetilde{H_0}(x)=(1-W^{-2})/2.$
The resulting recurrence equation is
\begin{eqnarray}
\label{genus0}
H^0_{(1^d)}=\frac1{d(d-1)}
\binom{2d-2}2
\sum_{i=1}^{d-1} i^2(d-i)^2
\binom{2d-4}{2i-2}
H^0_{(1^i)}H^0_{(1^{d-i})},
\end{eqnarray}
which is a well known recurrence found by Pandharipande (see \cite{v2}
Section 5.11 or \cite{fanp} p. 17).  Other (more complicated)
genus 0 recurrences can also be found in this manner.
}
\end{example}

%-----

\begin{example}[Closed form expressions for simple Hurwitz numbers]
{\rm
Closed form expressions for simple Hurwitz numbers
can be found for all genera (using the method of \cite{gj2} Cor.~4.1).
The expression for the genus~$g$ case can be obtained from Theorem~\ref{gjpf},
with the specializations of $p,s$ and $\phi_i$ given above,
and is the following. 

$$\frac{H_{(1^d)}^g}{(2d+2g-2)!}=\left[x^d\right]\widetilde{H_g}(x)
= \sum_{r=2g-1}^{5g-5} \sum_{n=r-1}^{r+g-1}
K_{n,g,r}\left( \left[x^d\right]\,\frac {w^n} { (1-w)^r }\right) $$
where
$$ K_{n,g,r}= \sum_{\substack{{\theta\models n} \\ {l(\theta)=r-2(g-1)}}}
(-1)^k\wi{\tau_{\theta_1}\tau_{\theta_2}\ldots\lambda_k}g $$
and $k=\sum_i(1-i)\theta_i+3g-3.$
Thus $K_{n,g,r}$ can be computed by Faber's program~\cite{fab}.
The remaining term is obtained by Lagrange inversion as
\begin{eqnarray*}
\left[x^d\right]\,\frac {w^n} { (1-w)^r } 
&=& \frac 1d \left[\mu^{d-1}\right]\,
\left(
\frac{n\mu^{n-1}}{(1-\mu)^r}
+\frac{r\mu^{n}}{(1-\mu)^{r+1}}
\right) e^{d\mu} \\
&=&
\sum_{i=0}^{d-n}\binom{r+i-1}{r-1}\frac{n\,d^{d-n-i-1}}{(d-n-i)!}+
\sum_{i=0}^{d-n-1}\binom{r+i}{r}\frac{r\,d^{d-n-i-2}}{(d-n-i-1)!}.
\end{eqnarray*}

For example,
for $\widetilde{H_3}(x),$ by Lagrange inversion,
\begin{eqnarray}
\frac{H^3_{(1^d)}}{(2d+4)!}
\nonumber
&=& {\frac {1}{1008}}\,{ A_4}(d) -{\frac {113}{10080}}\,{ A_5}(d)
+{\frac {2383 }{51840}}\,{ A_6}(d) -{\frac {16759}{181440}}\,{ A_7}(d) \\
&\mbox{}& +{\frac {227}{2304 }}\,{ A_8}(d)
-{\frac {557}{10368}}\,{ A_9}(d) +{\frac {245}{20736}}\,{ A_{10}}(d) \label{genus3}
\end{eqnarray}
where
$$A_k(d)=\frac kd \sum_{r=0}^{d-1}\binom{k+r}k\frac{d^{d-r-1}}{(d-r-1)!}.$$
This can be rewritten as
\begin{eqnarray*}
H^3_{(1^d)} &=& \frac {(2d+4)!}{2^53^3 9!}\,
\sum_{r=0}^{d-1}\frac{d^{d-r-2}}{(d-r-1)!} \binom{r+4}{5}(r+1)\\
&\mbox{}&\cdot\left (1225\,{r}^{4}+3770\,{r}^{3}+35\,{r}^{2}-2822\,r+1680\right ).
\end{eqnarray*}

It is clear that in general the simple Hurwitz numbers have the form
$$ H^g_{(1^d)}=(2d+2g-2)!\sum_{r=0}^{d-1}
\frac{d^{d-r-2}}{(d-r-1)!} P_g(d-r-1) $$
where $P_g(r)$ is a polynomial in $r$ of degree $5g-5.$

} \end{example}
%-------------------------------------------------------------------
% 
\noindent{\bf Acknowledgements.}  We would like to thank Carel Faber for
making available to us his {\sf Maple} program for evaluating Hodge integrals. 
We are grateful for conversations with
Gilberto Bini,
Robbert Dijkgraaf,
Tom Graber,
John Harer,
A.~Johan de Jong,
Rahul Pandharipande,
Malcolm Perry and
Wati Taylor.
DMJ would like to thank the Mathematics Departments at
Cambridge University, Duke University and MIT for their
hospitality during his sabbatical leave (1998/9) when much of this
work was carried out.

}  % end of parskip; it started just before the Introduction
%-------------------------------------------------------------------

\end{document}